\documentclass[9pt,a4paper]{article}
\usepackage{graphicx} \usepackage{array} \usepackage{subfigure} \usepackage{psfrag}
\usepackage{xy} \usepackage{amsmath} \usepackage{amssymb} \usepackage{fancyhdr}
\usepackage{setspace} \usepackage{rotating,epsfig} \usepackage{epsfig} 
\usepackage{verbatim} \usepackage{multicol} \usepackage{longtable}
\usepackage{color}\usepackage{times} \usepackage{epstopdf}
\usepackage{natbib} \usepackage{caption}

\def \BE {\begin{eqnarray*}} \def \EE {\end{eqnarray*}} % Allows &=& lineups in equations
\newcommand{\blue}[1]{\textcolor{blue}{#1}}
\newcommand{\red}[1]{\textcolor{red}{#1}}

\begin{document}
\pagenumbering{arabic}
\renewcommand{\thesection}{\arabic{section}}
\renewcommand\floatpagefraction{.9}
\renewcommand\topfraction{.9}
\renewcommand\bottomfraction{.9}
\renewcommand\textfraction{.1}

\makeatletter
\renewcommand{\section}{\@startsection{section}{1}{0mm}{-\baselineskip}{0.5\baselineskip}
{\normalfont\normalsize\scshape}}
\makeatother
\def \deg {$^{\circ}$}

% TITLE
\hrule
\begin{center}
\begin{large}\textbf{~~~~~~~~~~~~~~~Gaussian Processes for Regression:\newline
A Quick Introduction}\end{large}
\end{center}\begin{center}
~~~~~~~~~~~~~~~~~~~~~~~~~~M.~Ebden, August 2008\newline
Comments to mebden@gmail.com
\end{center}
\hrule

% THE BODY
\section{Motivation}

\begin{figure}[bH!]
 \begin{center}
  \includegraphics[height=7.5cm] {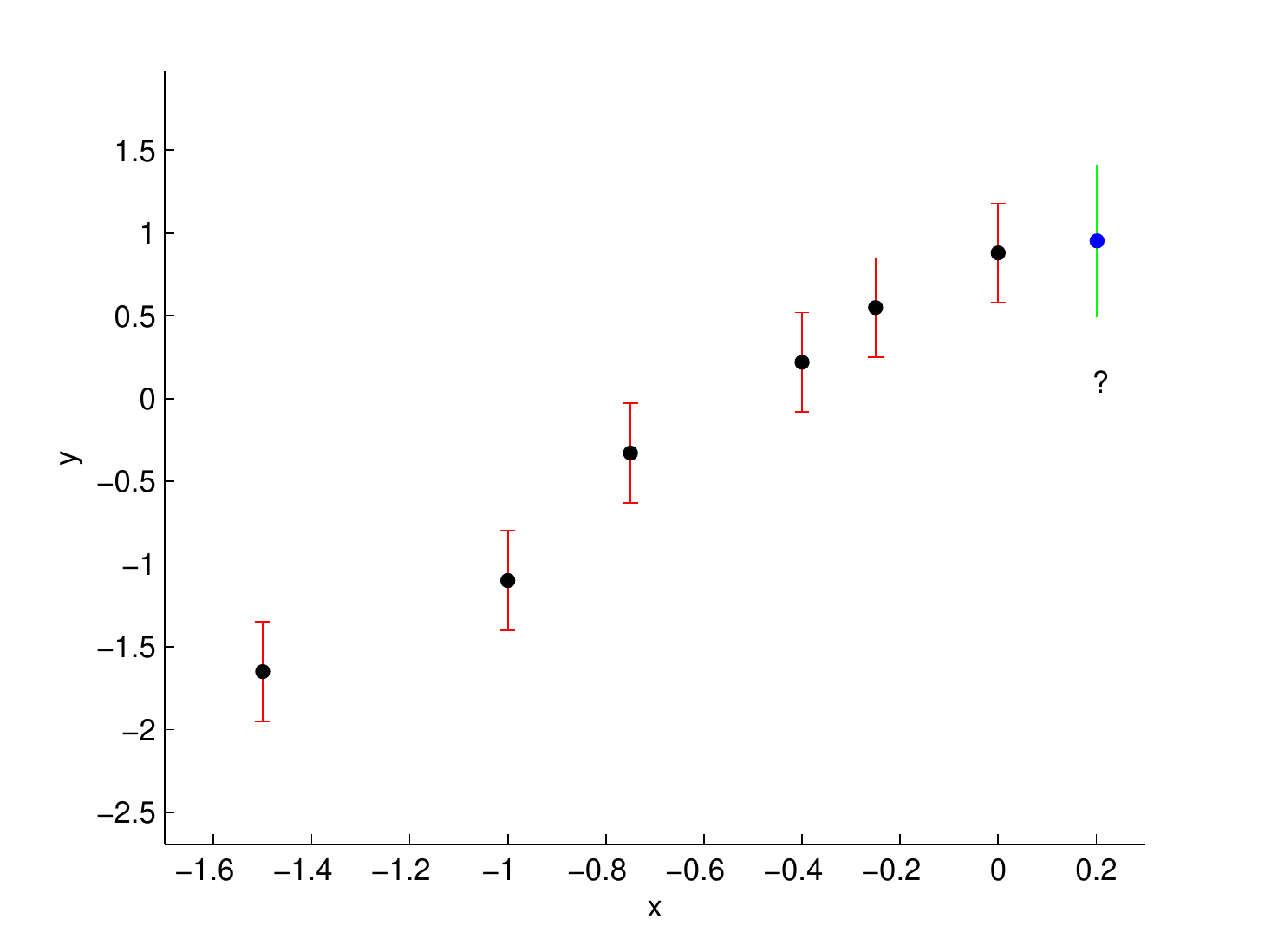} % needed for arxiv jp
  \caption{Given six noisy data points (error bars are indicated with vertical lines), we are interested in estimating a
seventh at $x_*=0.2$.}
  \label{fig:toyProblem}
 \end{center}
\end{figure}

Figure~\ref{fig:toyProblem} illustrates a typical example of a prediction problem: given some 
noisy observations of a dependent variable at certain values of the independent variable $x$, what
is our best estimate of the dependent variable at a new value, $x_*$?

If we expect the underlying function $f(x)$ to be linear, and can make some assumptions about the 
input data, we might use a least-squares method to fit a straight line (linear regression).  Moreover, if we suspect $f(x)$
may also be quadratic, cubic, or even nonpolynomial, we can use the principles of model selection to choose
among the various possibilities.

Gaussian process regression (GPR) is an even finer approach than this.
Rather than claiming $f(x)$ relates to some specific models (e.g.\ $f(x) = mx + c$), a Gaussian process can
represent $f(x)$ obliquely, but rigorously, by letting the data `speak' more clearly for themselves.  GPR is still a form
of supervised learning, but the training data are harnessed in a subtler way.

As such, GPR is a less `parametric' tool.  However, it's not completely free-form, and if we're unwilling to make
even basic assumptions about $f(x)$, then more general techniques should be considered, including those underpinned by the
principle of maximum entropy; Chapter~6 of \cite{H159} offers an introduction.

% GPtutorial.m
%training data (input and corresponding outputs) (x,y)
%new inputs, x*	 -> GPR -> new outputs
%arbitrary choice of covariance function k(x,x')

\section{Definition of a Gaussian Process}

Gaussian processes (GPs) extend multivariate Gaussian distributions to infinite dimensionality.
\blue{Formally, a Gaussian process generates data located throughout some domain such that any finite subset of
the range follows a multivariate Gaussian distribution.}  Now, the $n$ observations in an arbitrary
data set, $\mathbf{y} = \{y_1, \ldots, y_n\}$, can always be imagined as a single point sampled from some multivariate
($n$-variate) Gaussian distribution, after enough thought.  Hence, working backwards, this data set can be partnered with
a GP.  Thus GPs are as universal as they are simple.

Very often, it's assumed that the mean of this partner GP is zero everywhere.  What relates one
observation to another in such cases is just the \textit{covariance function}, $k(x,x')$.  A popular
choice is the `squared exponential',
\begin{equation}
k(x,x') = \sigma_f^2 \exp\bigg[\frac{-(x-x')^2}{2l^2}\bigg], \label{eq:se}
\end{equation}
where the maximum allowable covariance is defined as $\sigma_f^2$~--- this should be high for functions
which cover a broad range on the $y$ axis.  If $x \approx x'$, then $k(x,x')$ approaches
this maximum, meaning $f(x)$ is nearly perfectly correlated with $f(x')$.
This is good: for our function to look smooth, neighbours must be
alike.
Now if $x$ is distant from $x'$, we have instead $k(x,x') \approx 0$, i.e.\ the 
two points cannot `see' each other.  So, for example, during interpolation at new $x$ values, distant 
observations will have negligible effect.  How much effect this separation has will
depend on the length parameter, $l$, so there is much flexibility built into (\ref{eq:se}).

Not quite enough flexibility though: the data are often noisy as well, from measurement errors and so on.
Each observation $y$ can be thought of as related to an underlying function $f(x)$ through a Gaussian noise model:
\begin{equation}
 y = f(x) + \mathcal{N}(0,\sigma_n^2), \label{eq:noise}
\end{equation}
something which should look familiar to those who've done regression before.  Regression is the search for $f(x)$.
Purely for simplicity of exposition in the next page, we take the novel approach of folding the noise into
$k(x,x')$, by writing
\begin{equation}
k(x,x') = \sigma_f^2 \exp\bigg[\frac{-(x-x')^2}{2l^2}\bigg] + \sigma_n^2 \delta(x,x'), \label{eq:kxx}
\end{equation}
where $\delta(x,x')$ is the Kronecker delta function.  (When most people use Gaussian processes, they keep $\sigma_n$
separate from $k(x,x')$.  However, our redefinition of $k(x,x')$ is equally suitable for working with problems of the
sort posed in Figure~\ref{fig:toyProblem}.  So, given $n$ observations $\mathbf{y}$, our
objective is to predict $y_*$, not the `actual' $f_*$; their expected
values are identical according to (\ref{eq:noise}), but their variances differ owing to the observational
noise process.  e.g.\ in Figure~\ref{fig:toyProblem}, the expected value of $y_*$, and of $f_*$, is the dot at
$x_*$.)

To prepare for GPR, we calculate the covariance function, (\ref{eq:kxx}), among all
possible combinations of these points, summarizing our findings in three
matrices:
\begin{equation}
K = \begin{bmatrix}
k(x_1,x_1) & k(x_1,x_2) & \cdots & k(x_1,x_n) \\
k(x_2,x_1) & k(x_2,x_2) & \cdots & k(x_2,x_n) \\
\vdots & \vdots & \ddots & \vdots \\
k(x_n,x_1) & k(x_n,x_2) & \cdots & k(x_n,x_n)
\end{bmatrix} \label{eq:K}
\end{equation}
\begin{equation}
K_* = \begin{bmatrix}
k(x_*,x_1) & k(x_*,x_2) & \cdots & k(x_*,x_n)
\end{bmatrix}
~~~~~~~~~~~~~~ K_{**} = k(x_*,x_*). \label{eq:Ks}
\end{equation}
Confirm for yourself that the diagonal elements of $K$ are $\sigma_f^2 + \sigma_n^2$, and that
its extreme off-diagonal elements tend to zero when $x$ spans a large
enough domain.

\section{How to Regress using Gaussian Processes}

\begin{figure}
 \begin{center}
  \includegraphics[height=7.6cm] {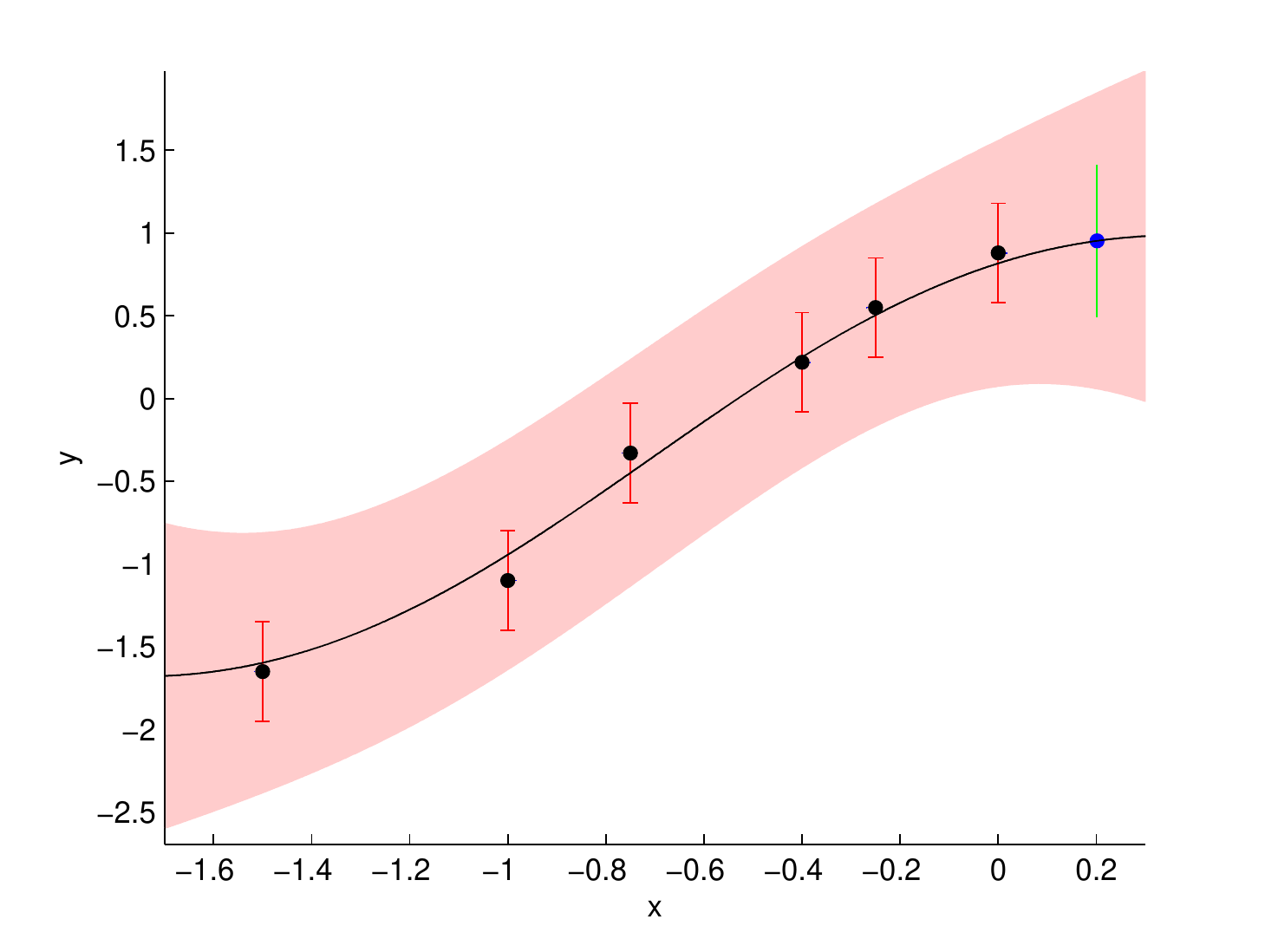}
  \caption{The solid line indicates an estimation of $y_*$ for 1,000 values of $x_*$.  Pointwise 95\%
confidence intervals are shaded.}
  \label{fig:prettier}
 \end{center}
\end{figure}

Since the key assumption in GP modelling is that our data can be represented
as a sample from a multivariate Gaussian distribution, we have that
\begin{equation}
\begin{bmatrix}
\mathbf{y} \\
y_* \end{bmatrix} \sim \mathcal{N} \bigg( \mathbf{0},
\begin{bmatrix}
K & K_*^\mathrm{T} \\
K_* & K_{**}
\end{bmatrix}
\bigg), \label{eq:kmat}
\end{equation}
where $\mathrm{T}$ indicates matrix transposition.  We are of course interested in the conditional probability 
$p(y_*|\mathbf{y})$: ``given the data, how likely is a certain prediction for $y_*$?''.  As explained more slowly in the Appendix,
the probability follows a Gaussian distribution:
\begin{equation}
y_* | \mathbf{y} \sim \mathcal{N} (K_* K^{-1} \mathbf{y}, ~K_{**} - K_* K^{-1} K_*^\mathrm{T}). \label{eq:ycond}
\end{equation}
Our best estimate for $y_*$ is the mean of this distribution:
\begin{equation}
\overline{y}_* = K_* K^{-1} \mathbf{y}, \label{eq:ybar}
\end{equation}
and the uncertainty in our estimate is captured in its variance:
\begin{equation}
\mathrm{var}(y_*) = K_{**} - K_* K^{-1} K_*^\mathrm{T}. \label{eq:var_y}
\end{equation}
We're now ready to tackle the data in Figure~\ref{fig:toyProblem}.
\begin{enumerate}
\item There are $n=6$ observations $\mathbf{y}$, at
\[ \mathbf{x} = 
\begin{bmatrix}
-1.50 & -1.00 & -0.75 & -0.40 & -0.25 & 0.00
\end{bmatrix}
. \]
We know $\sigma_n = 0.3$ from the error bars.  With judicious choices of $\sigma_f$ and $l$
(more on this later), we have enough to calculate a covariance matrix using (\ref{eq:K}):
\[K =
\begin{bmatrix}
\red{1.70} & 1.42 & 1.21 & 0.87 & 0.72 & 0.51 \\ 
1.42 & \red{1.70} & 1.56 & 1.34 & 1.21 & 0.97 \\ 
1.21 & 1.56 & \red{1.70} & 1.51 & 1.42 & 1.21 \\ 
0.87 & 1.34 & 1.51 & \red{1.70} & 1.59 & 1.48 \\ 
0.72 & 1.21 & 1.42 & 1.59 & \red{1.70} & 1.56 \\ 
0.51 & 0.97 & 1.21 & 1.48 & 1.56 & \red{1.70}
\end{bmatrix}.
\] % GPtutorial.m
From (\ref{eq:Ks}) we also have $K_{**} = 1.70$ and
\[K_* = 
\begin{bmatrix}
%0.31 & 0.68 & 0.92 & 1.25 & 1.38 & 1.54  -- old, wrong values.  see email of 23/4/12
0.38 & 0.79 & 1.03 & 1.35 & 1.46 & 1.58
\end{bmatrix}.
\]

\item From (\ref{eq:ybar}) and (\ref{eq:var_y}), $\overline{y}_* = 0.95$ and
$\mathrm{var}(y_*) = 0.21$.

\item Figure~\ref{fig:toyProblem} shows a data point with a question mark
underneath, representing
the estimation of the dependent variable at $x_*=0.2$.
\end{enumerate}
We can repeat the above procedure for various other points spread over some portion of the $x$ axis, as shown in
Figure~\ref{fig:prettier}.  (In fact, equivalently, we could avoid the
repetition by performing the above procedure once
with suitably larger $K_*$ and $K_{**}$ matrices.  In this case, since there are 1,000 test points spread over the $x$
axis, $K_{**}$ would be of size 1,000 $\times$ 1,000.)  Rather than plotting simple error bars, we've decided to plot
$\overline{y}_* \pm 1.96 \sqrt{\mathrm{var}(y_*)}$, giving a 95\% confidence interval.

% GPtutorial.m

\section{GPR in the Real World}

The reliability of our regression is dependent on how well we select the
covariance function.
Clearly if its parameters~--- call them $\boldsymbol{\theta} = \{l,
\sigma_f,
\sigma_n\}$~--- are 
not chosen sensibly, the result is nonsense.  Our maximum \textit{a posteriori} estimate of
$\boldsymbol{\theta}$ occurs when $p(\boldsymbol{\theta}|\mathbf{x},\mathbf{y})$ is
at its greatest.  Bayes' theorem tells us that, assuming we have little
prior knowledge about what $\boldsymbol{\theta}$
should be, this corresponds to maximizing $\log
p(\mathbf{y}|\mathbf{x},\boldsymbol{\theta})$, given by
\begin{equation}
\log p(\mathbf{y}|\mathbf{x},\boldsymbol{\theta}) = -\frac{1}{2}
\mathbf{y}^\mathrm{T} K^{-1} \mathbf{y} - \frac{1}{2} \log |K| - \frac{n}{2} \log 2 \pi. \label{eq:ml}
\end{equation}
Simply run your favourite multivariate optimization algorithm (e.g.\ conjugate gradients,
Nelder-Mead simplex, etc.) on this equation and you've found a pretty good choice for
$\boldsymbol{\theta}$; in our example, $l=1$ and $\sigma_f = 1.27$.

It's only ``pretty good'' because, of course, Thomas Bayes is rolling in 
his grave.  Why commend just one answer for $\boldsymbol{\theta}$, when
you can integrate everything over the many different possible choices for 
$\boldsymbol{\theta}$?  Chapter~5 of \citet{H152} presents the
equations necessary in this case.

\begin{figure}
 \begin{center}
  \includegraphics[width=14.45cm] {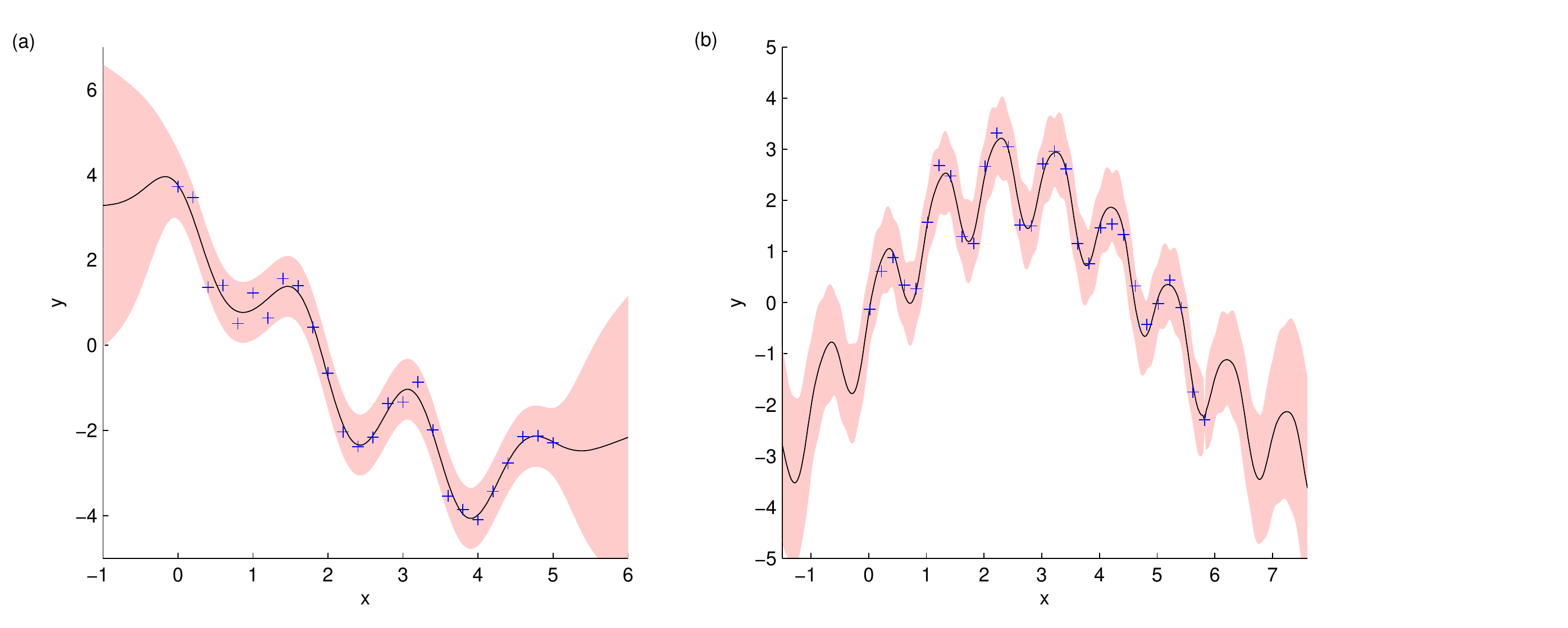}
  \caption{~Estimation of $y_*$ (solid line) for a function with \textbf{(a)} short-term and long-term dynamics, and
\textbf{(b)}~long-term dynamics and a periodic element.  Observations are shown as crosses.}
  \label{fig:more_fcns}
 \end{center}
\end{figure}

Finally, if you feel you've grasped the toy problem in
Figure~\ref{fig:prettier}, the next two
examples handle more complicated cases.  Figure~\ref{fig:more_fcns}(a), in addition to
a long-term downward trend, has some fluctuations, so we might use a more 
sophisticated covariance function:
\begin{equation}
k(x,x') = {\sigma_f}_1^2 \exp\bigg[\frac{-(x-x')^2}{2l_1^2}\bigg] + {\sigma_f}_2^2
\exp\bigg[\frac{-(x-x')^2}{2l_2^2}\bigg] + \sigma_n^2 \delta(x,x').
\end{equation}
The first term takes into account the small vicissitudes of the dependent
variable, and
the second term has a longer length parameter ($l_2 \approx 6 l_1$) to
represent its long-term
trend.  Covariance functions can be grown in this way \textit{ad
infinitum}, to suit the complexity of your particular data.

% GPtutorial2.m, GPtutorial3.m

The function looks as if it might contain a periodic element, but it's difficult to
be sure.  Let's consider another function, which we're told has a
periodic element.  The solid line in Figure~\ref{fig:more_fcns}(b) was regressed with the following covariance
function:
\begin{equation}
k(x,x') = \sigma_f^2 \exp\bigg[\frac{-(x-x')^2}{2l^2}\bigg] + \exp\{-2 \sin^2[\nu \pi
(x-x')]\} + \sigma_n^2 \delta(x,x'). 
\end{equation}
The first term represents the hill-like trend over the long term, and
the second term gives periodicity with frequency $\nu$.  This is the first time we've encountered a
case where $x$ and $x'$ can be distant and yet still `see' each other (that is, $k(x,x') \not\approx 0$ for
$x \gg x'$).

What if the dependent variable has other dynamics which, \textit{a priori}, you 
expect to appear?  There's no limit to how complicated $k(x,x')$ can be, provided $K$ is positive
definite.  Chapter 4 of \citet{H152} offers a good outline of the range of covariance functions you should keep in
your toolkit.

``Hang on a minute,'' you ask, ``isn't choosing a covariance function from a toolkit a lot like
choosing a model type, such as linear versus cubic~--- which we discussed at the outset?''  Well,
there are indeed similarities.  In fact, there is no way to perform regression without imposing at
least a modicum of structure on the data set; such is the nature of generative modelling.
However, it's worth repeating that Gaussian processes do allow the data to speak very clearly.
For example, there exists excellent theoretical justification for the use of (\ref{eq:se}) in many settings (\citet{H152},
Section~4.3).  You will still want to investigate carefully which covariance functions are appropriate for your data set.
Essentially, choosing among alternative functions is a way of reflecting various forms of prior
knowledge about the \textit{physical process} under investigation.

\section*{}
\section*{}
\section{Discussion}

We've presented a brief outline of the mathematics of GPR, but practical implementation of the above
ideas requires the solution of a few algorithmic hurdles as opposed to those of data analysis.  If you aren't a good computer
programmer, then the code for Figures~\ref{fig:toyProblem} and~\ref{fig:prettier} is at
\textit{github.com/mebden/GPtutorial}, and more general code can be found
at \textit{www.gaussianprocess.org/gpml}.
%\textit{www.robots.ox.ac.uk/}$\sim$\textit{mebden/reports/GPtutorial.zip}

We've merely scratched the surface of a powerful technique \citep{H160}. 
First, although the focus has been on one-dimensional inputs, it's simple
to accept those of higher dimension.  Whereas $x$ would then change from a
scalar to a vector, $k(x,x')$ would remain a scalar and so the maths
overall would be virtually unchanged.  Second, the zero vector
representing the mean of the multivariate Gaussian distribution in
(\ref{eq:kmat}) can be replaced with functions of $x$.  Third, in addition
to their use in regression, GPs are applicable to integration, global
optimization, mixture-of-experts models, unsupervised learning models, and
more~--- see Chapter~9 of \citet{H152}.  The next tutorial will focus on
their use in \textit{classification}.

\bibliographystyle{Chicago}
\bibliography{GPtute}

\section*{Appendix}

Imagine a data sample $\mathbf{d}$ taken from some multivariate Gaussian distribution with zero mean and a covariance
given by matrix $D$.  Now decompose $\mathbf{d}$ arbitrarily into two consecutive subvectors $\mathbf{a}$ and $\mathbf{b}$~--- in
other words, writing $\mathbf{d} \sim \mathcal{N} ( \mathbf{0}, D)$ would be the same as writing
\begin{equation}
\begin{bmatrix}
\mathbf{a} \\
\mathbf{b}
\end{bmatrix} \sim \mathcal{N} \bigg( \mathbf{0},
\begin{bmatrix}
A & C^\mathrm{T} \\
C & B
\end{bmatrix}
\bigg), \label{eq:genmat}
\end{equation}
where $A$, $B$, and $C$ are the corresponding bits and pieces that make up $D$.

Interestingly, the conditional distribution of $\mathbf{b}$ given $\mathbf{a}$ is itself Gaussian-distributed.  If the
covariance matrix $D$ were diagonal or even block diagonal, then knowing $\mathbf{a}$ wouldn't tell us anything about
$\mathbf{b}$: specifically, $\mathbf{b} | \mathbf{a} \sim
\mathcal{N} (\mathbf{0}, B)$.  On the other hand, if $C$ were nonzero, then some matrix algebra leads us to
\begin{equation}
\mathbf{b} | \mathbf{a} \sim \mathcal{N} (CA^{-1}\mathbf{a}, ~B - C A^{-1} C^\mathrm{T} ).\label{eq:ba}
\end{equation}
The mean, $CA^{-1}\mathbf{a}$, is known as the `matrix of regression coefficients', and the variance, $B - C A^{-1}
C^\mathrm{T}$, is the `Schur complement of $A$ in $D$'.

In summary, if we know some of $\mathbf{d}$, we can use that to inform our estimate of what the rest of $\mathbf{d}$ might be,
thanks to the revealing off-diagonal elements of $D$.

\newpage

\hrule
\begin{center}
\begin{large}\textbf{~~~~~~~~~~~~~Gaussian Processes for Classification:\newline
A Quick Introduction}\end{large}
\end{center}\begin{center}
~~~~~~~~~~~~~~~~~~~~~~~~~~M.~Ebden, August 2008\newline
Prerequisite reading: \textit{Gaussian Processes for Regression}
\end{center}
\hrule

% THE BODY
\setcounter{section}{0} \setcounter{equation}{0} \setcounter{figure}{0}
\setcounter{footnote}{0} \setcounter{table}{0}
\setcounter{page}{1}

\section{Overview}

As mentioned in the previous document, GPs can be applied to problems other than regression.  For example, if the output of a GP is squashed onto the range $[0,1]$, it can represent the
\textit{probability} of a data point belonging to one of say two types, and voil\`{a}, we can ascertain classifications.  This is the
subject of the current document.

The big difference between GPR and GPC is how the output data, $\mathbf{y}$, are linked to the underlying function
outputs, $\mathbf{f}$.  They are no longer connected simply via a noise process as in (\ref{eq:noise}) in the previous document, but
are instead now discrete: say $y=1$ precisely for one class and $y=-1$ for the other. In
principle, we could try fitting a GP that produces an output of approximately $1$ for some values of $x$ and approximately $-1$ for
others, simulating this discretization.  Instead, we interpose the GP between the data and a
squashing function; then, classification of a new data point $x_*$ involves two steps instead of one:
\begin{enumerate}
\item Evaluate a `latent function' $f$ which models qualitatively how the likelihood of one class versus the other
changes over the $x$ axis.  This is the GP.
\item Squash the output of this latent function onto $[0,1]$ using any sigmoidal function, $\pi(f) = \text{prob}(y=1|f)$.
\end{enumerate}
Writing these two steps schematically,
\begin{center}
\framebox{data, $x_*$} $\xrightarrow{\text{GP}}$ \framebox{latent function, $f_*|x_*$}
$\xrightarrow{\text{sigmoid}}$ \framebox{class probability, $\pi(f_*)$}.
\end{center}
The next section will walk you through more slowly how such a classifier operates.  Section~\ref{sec:trainGPC} explains
how to train the classifier, so perhaps we're presenting things in reverse order!  Section~\ref{sec:multiC} handles
classification when there are more than two classes.

Before we get started, a quick note on $\pi(f)$.  Although other forms will do, here we will prescribe it to be the cumulative Gaussian
distribution, $\Phi(f)$.  This $S$-shaped function satisfies our needs, mapping high $f$ 
into $\pi(f) \approx 1$, and low $f$ into $\pi(f) \approx 0$.

A second quick note, revisiting (\ref{eq:kmat}) and (\ref{eq:ycond}) in the first document: confirm for yourself that, if there were
no noise ($\sigma_n = 0$), the two equations could be rewritten as
\begin{equation}
\begin{bmatrix}
\mathbf{f} \\
f_* \end{bmatrix} \sim \mathcal{N} \bigg( \mathbf{0},
\begin{bmatrix}
K & K_*^\mathrm{T} \\
K_* & K_{**}
\end{bmatrix}
\bigg) \label{eq:kmatf}
\end{equation}
and
\begin{equation}
f_* | \mathbf{f} \sim \mathcal{N} (K_* K^{-1} \mathbf{f}, ~K_{**} - K_* K^{-1} K_*^\mathrm{T}). \label{eq:fcond}
\end{equation}
\begin{comment}
\begin{equation}
\overline{f}_* = K_* K^{-1} \mathbf{f}, \label{eq:fbar}
\end{equation}
and 
\begin{equation}
\mathrm{var}(f_*) = K_{**} - K_* K^{-1} K_*^\mathrm{T}. \label{eq:varf}
\end{equation}
\begin{equation}
\overline{f}_* = K_* K^{-1} \hat{\mathbf{f}}, \label{eq:fbarC}
\end{equation}
and the variance works out as
\begin{equation}
\mathrm{var}(f_*) = K_{**} - K_* (K')^{-1} K_*^\mathrm{T}. \label{eq:varfC}
\end{equation}
\end{comment}
\[ \]
\section*{}

\section{Using the Classifier}

Suppose we've trained a classifier from $n$ input data, $\mathbf{x}$, and their corresponding expert-labelled output data, $\mathbf{y}$.
And suppose that in the process we formed some GP outputs $\mathbf{f}$ corresponding to these data, which have some uncertainty
but mean values given by
$\hat{\mathbf{f}}$. We're now ready to input a new data point, $x_*$, in the left side of our schematic, in order to determine at the other
end the probability $\pi_*$ of its class membership.

\begin{figure}
 \begin{center}
  \includegraphics[width=12cm] {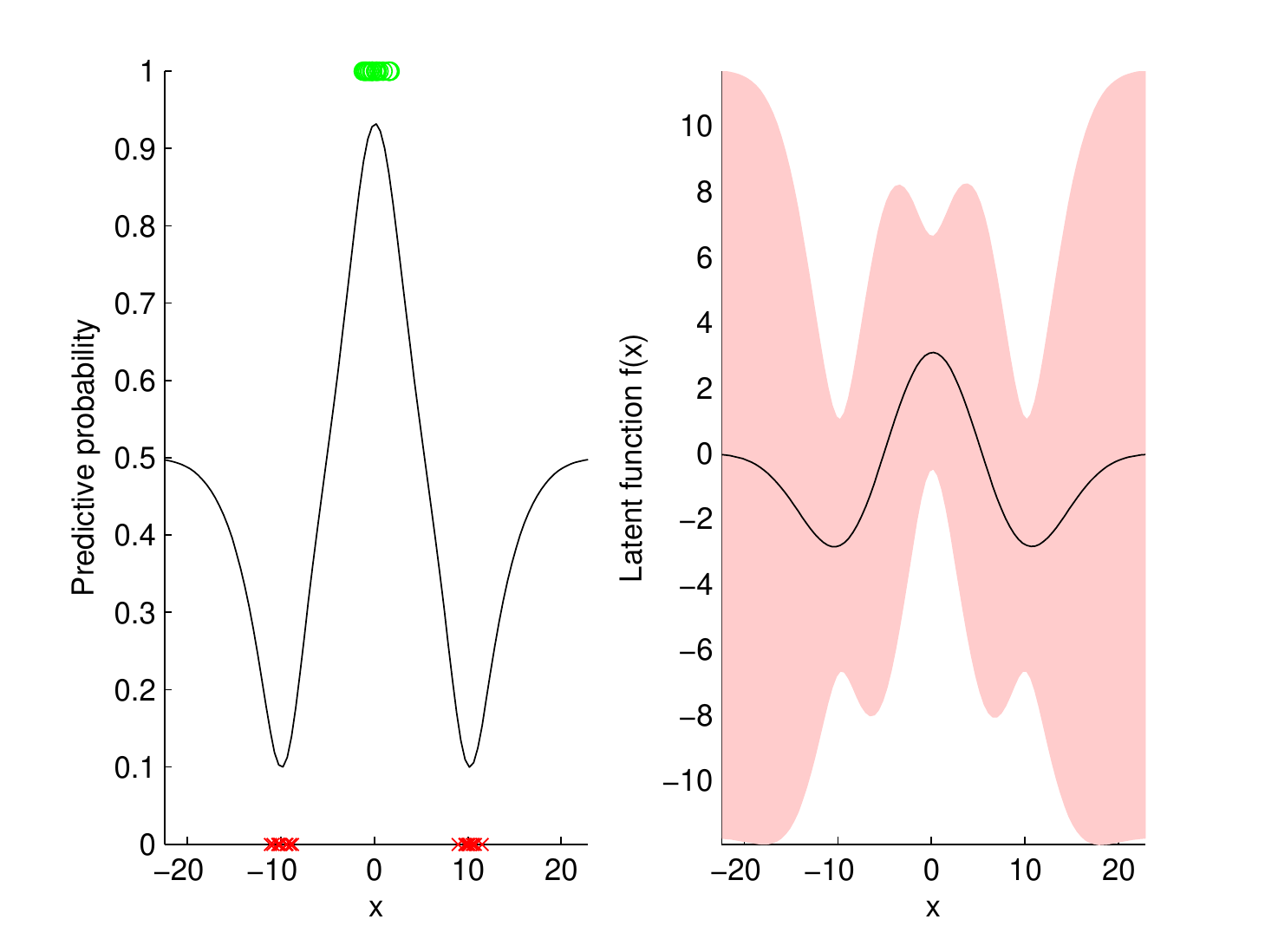}
  \caption{\textbf{(a)}~Toy classification dataset, where circles and crosses indicate class membership of the input (training)
data $\mathbf{y}$, at locations $\mathbf{x}$.  $\mathbf{y}$ is $1$ for one class and $-1$ for another, but for illustrative purposes we
pretend $y=0$ instead of $-1$ in this figure.  The solid line is the (mean) probability $\overline{\pi}_* =
\text{prob}(y_*=1|x_*)$, i.e.\ the `answer' to our problem after successfully performing GPC.  \textbf{(b)}~The
corresponding distribution of the latent function $f$, not constrained to lie between 0 and 1.}
  \label{fig:gpc}
 \end{center}
\end{figure}

In the first
step, finding the probability $p(f_*|\mathbf{f})$ is similar to GPR, i.e.\ we adapt (\ref{eq:fcond}):
\begin{equation}
p(f_* | \mathbf{f}) = \mathcal{N} (K_* K^{-1} \hat{\mathbf{f}}, ~K_{**} - K_* (K')^{-1} K_*^\mathrm{T}). \label{eq:fcond2}
\end{equation}
($K'$ will be explained soon, but for now consider it to be very similar to $K$.)  In the second step, we
squash $f_*$ to find the probability of class membership, $\pi_* = \pi(f_*) = \Phi(f_*)$.  The expected
value is
\begin{equation}
\overline{\pi}_* = \int \pi(f_*) p(f_* | \mathbf{f}) df_*.
\end{equation}
This is the integral of a cumulative Gaussian times a Gaussian, which can be solved analytically.  By Section~3.9 of
\cite{H152}, the solution is:
\begin{equation}
\overline{\pi}_* = \Phi \bigg( \frac{\overline{f}_*}{\sqrt{1 + \mathrm{var}(f_*)}} \bigg).
\end{equation}
An example is depicted in Figure~\ref{fig:gpc}.

% GPtutorial?.m, GPtutorial?.m

\section{Training the GP in the Classifier}\label{sec:trainGPC}

Our objective now is to find $\hat{\mathbf{f}}$ and $K'$, so that we know
everything about the GP producing (\ref{eq:fcond2}), the first step of the
classifier.  The second step of the classifier does not require training as it's a fixed sigmoidal function.

Among the many GPs which could be partnered with our data set, naturally we'd like to compare
their usefulness quantitatively.  Considering the outputs $\mathbf{f}$ of a certain GP, how likely they are to be appropriate for
the training data can be decomposed using Bayes' theorem:
\begin{equation}
p(\mathbf{f}|\mathbf{x},\mathbf{y}) = \frac{p(\mathbf{y}|\mathbf{f})p(\mathbf{f}|\mathbf{x})}{p(\mathbf{y}|\mathbf{x})}.
\label{eq:bayes2}
\end{equation}
Let's focus on the two factors in the numerator.  Assuming the data set is i.i.d., 
\begin{equation}
p(\mathbf{y}|\mathbf{f}) = \prod_{i=1}^n p(y_i|f_i). \end{equation}
Dropping the subscripts in the product, $p(y|f)$ is informed by our 
sigmoid function, $\pi(f)$.  Specifically, 
$p(y=1|f)$ is $\pi(f)$ by
definition, and to complete the picture, $p(y=-1|f) = 1 - \pi(f)$.  A
terse way of combining these two cases is to write $p(y|f) = \Phi(yf)$.  

The second factor in the numerator is $p(\mathbf{f}|\mathbf{x})$. This is related to the output
of the first step of our schematic drawing, but first we're interested in the value of $p(\mathbf{f}|\mathbf{x})$ which maximizes
the posterior
probability $p(\mathbf{f}|\mathbf{x},\mathbf{y})$.  This occurs when the derivative of (\ref{eq:bayes2}) with respect to $\mathbf{f}$ is zero, or
equivalently and more simply, when the derivative of its logarithm is zero.  Doing this, and using the same logic that
produced (\ref{eq:ml}) in the previous document, we find that
\begin{equation}
\hat{\mathbf{f}} = K \nabla \log p(\mathbf{y}|\hat{\mathbf{f}}), \label{eq:fhat}
\end{equation}
where $\hat{\mathbf{f}}$ is the best $\mathbf{f}$ for our problem.  Unfortunately, $\hat{\mathbf{f}}$ appears on both sides of
the equation, so we make an initial guess (zero is fine) and go through a few iterations.  The answer to (\ref{eq:fhat}) can be used
directly in (\ref{eq:fcond2}), so we've found one of the two quantities we seek therein.
%\begin{equation}
%\log p(\mathbf{f}|\mathbf{x}) = -\frac{1}{2}\mathbf{f}^\mathrm{T} K^{-1} \mathbf{f} - \frac{1}{2} \log |K| - \frac{n}{2} \log 2 \pi. \label{eq:ml2}
%\end{equation}

The variance of $\mathbf{f}$ is given by the negative second derivative of the logarithm of 
(\ref{eq:bayes2}), which turns out to be $(K^{-1}+W)^{-1}$, with $W = - \nabla \nabla \log p(\mathbf{y}|\mathbf{f})$.  Making
a \textit{Laplace approximation}, we pretend $p(\mathbf{f}|\mathbf{x},\mathbf{y})$ is Gaussian distributed, i.e.
\begin{equation}
p(\mathbf{f}|\mathbf{x},\mathbf{y}) \sim q(\mathbf{f}|\mathbf{x},\mathbf{y}) = \mathcal{N} (\hat{\mathbf{f}},
(K^{-1}+W)^{-1}).
\end{equation}
(This assumption is occasionally inaccurate, so if it yields poor classifications, better ways of
characterizing the uncertainty in $\mathbf{f}$ should be considered, for example via expectation propagation.)

Now for a subtle point.  The fact that $\mathbf{f}$ can vary means that using (\ref{eq:fcond}) directly is inappropriate:
in particular,
its mean is correct but its variance no longer tells the whole story.  This is why we use the adapted version, (\ref{eq:fcond2}),
with $K'$ instead of $K$.  Since the varying quantity in (\ref{eq:fcond}), $\mathbf{f}$, is being multiplied by $K_*K^{-1}$, we add $K_*K^{-1} \text{cov}(\mathbf{f})
K^{-1}K_*^{\mathrm{T}}$ to the variance in (\ref{eq:fcond}).  Simplification leads to (\ref{eq:fcond2}), in which $K' = K
+ W^{-1}$.

With the GP now completely specified, we're ready to use the classifier as described in the previous section.
\newline
\newline
\noindent{\textbf{GPC in the Real World}}
\newline
As with GPR, the reliability of our classification is dependent on how well we select the covariance function
in the GP in our first step.  The parameters are $\boldsymbol{\theta} = \{l, \sigma_f \}$, one fewer now because $\sigma_n = 0$.
However, as usual, $\boldsymbol{\theta}$ is optimized by maximizing $p(\mathbf{y}|\mathbf{x},\boldsymbol{\theta})$, or (omitting
$\boldsymbol{\theta}$ on the righthand side of the equation),
\begin{equation}
p(\mathbf{y}|\mathbf{x},\boldsymbol{\theta}) = \int p(\mathbf{y}|\mathbf{f}) p(\mathbf{f}|\mathbf{x}) d\mathbf{f}.
\end{equation}
This can be simplified, using a Laplace approximation, to yield
\begin{equation}
p(\mathbf{y}|\mathbf{x},\boldsymbol{\theta}) = -\frac{1}{2}\hat{\mathbf{f}}^\mathrm{T}K^{-1}\hat{\mathbf{f}} + \log
p(\mathbf{y}|\hat{\mathbf{f}}) - \frac{1}{2} \log (|K|\cdot|K^{-1} + W|). \label{eq:ml2}
\end{equation}
This is the equation to run your favourite optimizer on, as performed in GPR.

%The only other difference is that the objective function passed to our optimizer will take longer
%to run than before, as it must iterate (\ref{eq:fhat}) before it can calculate $p(\mathbf{y}|\mathbf{x}$).

\section{Multi-Class GPC}\label{sec:multiC}

We've described binary classification, where the number of possible classes, $C$, is just two.  In the case of $C>2$ classes, 
one approach is to fit an $f$ for each class.  In the first of the two steps of classification, our GP values are concatenated as
\begin{equation}
\mathbf{f} = (f_1^1, \ldots, f_n^1, f_1^2, \ldots, f_n^2, \ldots, f_1^C, \ldots, f_n^C)^\mathrm{T}.
\end{equation}
Let $\mathbf{y}$ be a vector of the same length as $\mathbf{f}$ which, for each $i=1,\ldots,n$, is $1$ for the class which is
the label and $0$ for the other $C-1$ entries.  Let $K$ grow to being block diagonal in the matrices $K^1, \ldots, K^C$.  So the first
change we see for $C>2$ is a lengthening of the GP.  Section~3.5 of \cite{H152} offers hints on keeping the computations manageable.

The second change is that the (merely one-dimensional) cumulative Gaussian distribution is no longer sufficient to describe the squashing
function in our classifier; instead we use the softmax function.  For the $i$th data point,
\begin{equation}
p(y_i^c|\mathbf{f}_i) = \pi_i^c = \frac{\exp (f_i^c)}{\sum_{c'} \exp (f_i^{c'}) }
\end{equation}
where $\mathbf{f}_i$ is a nonconsecutive subset of $\mathbf{f}$, viz.\ $\mathbf{f}_i = \{f_i^1,f_i^2,\ldots,f_i^C \}$.  We can summarize our results
with $\boldsymbol{\pi} = \{\pi_1^1, \ldots, \pi_n^1, \pi_1^2, \ldots, \pi_n^2, \ldots, \pi_1^C, \ldots, \pi_n^C \}$.

Now that we've presented the two big changes needed to go from binary- to multi-class GPC, we continue as before.  Setting to zero
the derivative of the logarithm of the components in (\ref{eq:bayes2}), we replace (\ref{eq:fhat}) with
\begin{equation}
\hat{\mathbf{f}} = K (\mathbf{y} - \hat{\boldsymbol{\pi}}). \label{eq:fhatMulti}
\end{equation}
The corresponding variance is $(K^{-1}+W)^{-1}$ as before, but now $W = \text{diag}(\boldsymbol{\pi}) - \Pi \Pi^\mathrm{T}$, where $\Pi$ is
a $Cn \times n$ matrix obtained by stacking vertically the diagonal matrices $\text{diag}(\boldsymbol{\pi}^c)$, if $\boldsymbol{\pi}^c$ is the
subvector of $\boldsymbol{\pi}$ pertaining to class $c$.

With these quantities estimated, we have enough to generalize (\ref{eq:fcond2}) to
\begin{equation}
p(f_*^c | \mathbf{f}) = \mathcal{N} \big(K_*^c (K^c)^{-1} \hat{\mathbf{f}}^c, ~\text{diag}(K_{**}) - (K_*^c)^\mathrm{T} (K^c + (W^c)^{-1})^{-1}
(K_*^c)^\mathrm{T} \big), \label{eq:fcond3}
\end{equation}
where $f_*^c$, $K_*^c$, and $W^c$ represent the class-relevant information only.  Finally, (\ref{eq:ml2}) is replaced with
\begin{equation}
p(\mathbf{y}|\mathbf{x},\boldsymbol{\theta}) = -\frac{1}{2}\hat{\mathbf{f}}^\mathrm{T} K^{-1} \hat{\mathbf{f}} + \mathbf{y}^\mathrm{T}
\hat{\mathbf{f}} - \sum_{i=1}^n \log \Bigg[ \sum_{c=1}^C \exp \hat{f}_i^c \Bigg] - \frac{1}{2} \log \big(|K|\cdot|K^{-1} + W|\big). \label{eq:ml3}
\end{equation}
We won't present an example of multi-class GPC, but hopefully you get the idea.

\section{Discussion}

As with GPR, classification can be extended to accept $x$ values with multiple dimensions, while keeping most of the mathematics unchanged.  
Other possible extensions include using the expectation propagation method in lieu of the Laplace approximation as mentioned previously,
putting confidence intervals on the classification probabilities, calculating the derivatives of (\ref{eq:ml3}) to aid the optimizer, or
using the variational Gaussian process classifiers described in \citet{H160}, to name but four extensions.

Second, we repeat the Bayesian call from the previous document to integrate over a range of possible covariance function
parameters.  This should be done regardless of how much prior knowledge is available~--- see for example Chapter~5 of \cite{H159} on how to choose priors in
the most opaque situations.

Third, we've again spared you a few practical algorithmic details; computer code is available at \textit{www.gaussianprocess.org/gpml},
with examples.

\section*{Acknowledgments}

Thanks are due to Prof Stephen Roberts and members of the Pattern Analysis Research Group, as well as the ALADDIN project (www.aladdinproject.org).

\bibliography{GPtute}

\newpage
\setcounter{page}{1}
% TITLE
\hrule
\begin{center}
\begin{large}\textbf{Gaussian Processes for Dimensionality Reduction:\\
A Quick Introduction}\end{large}
\end{center}\begin{center}
M.~Ebden, August 2015\\
Prerequisite reading: \textit{Gaussian processes for regression}
\end{center}
\hrule
\setcounter{section}{0} \setcounter{equation}{0} \setcounter{figure}{0}
\setcounter{footnote}{0} \setcounter{table}{0}

% THE BODY

\section*{}

Suppose you wanted to learn a low-dimensional representation of a dataset for ease of interpretation, sort of like how your shadow is a simplified, two-dimensional representation of a three-dimensional object. Let the original data be matrix $\mathbf{Y}$ with $n$ rows (observations) of $d$ dimensions each, and let the new representation be $\mathbf{X}$ with $n$ rows of $q$ dimensions each ($q < d$). 
%= [\mathbf{y}_1, \ldots , \mathbf{y}_n]^\text{T}
%first assume a model of how to go from $\mathbf{X}$ to $\mathbf{Y}$. Specifically we'll

To learn this low-dimensional representation, we'll assume that for the $i$th dimension of $\mathbf{Y}$ the $n$ elements in $\mathbf{y}_{:,i}$ are a sample from a Gaussian process based on the low-dimensional space. Specifically, we'll use a zero-mean Gaussian process, $\mathcal{GP} \left( \mathbf{0}_{n \times 1}, k(\mathbf{x},\mathbf{x}') \right)$, keeping the same model for each dimension of $\mathbf{Y}$. We'll choose for $k(\cdot)$ the squared exponential kernel (a.ka.\ RBF kernel) in order to ensure that points close in $\mathbf{X}$ will be close in $\mathbf{Y}$. Renaming the noise variance from $\sigma_n^2$ (in the first report) to $1/\beta$, the kernel is:
\[k(\mathbf{x},\mathbf{x}') = \sigma^2 \exp\bigg[\frac{-|\mathbf{x}-\mathbf{x}'|^2}{2l^2}\bigg] + \frac{\delta(x,x')}{\beta}. \label{eq:kxx2}
\]
The covariance matrix $\mathbf{K}$ for this GP is constructed as per (4) in the first report, and we have no need for a $\mathbf{K}_*$ or $\mathbf{K}_{**}$ from (5) because there are no points in $\mathbf{Y}$ to predict. The tuple of kernel parameters this time is $\boldsymbol{\theta} = \{\sigma,l,\beta\}$.

We'll further assume that the GPs behind each of the $d$ dimensions have been sampled independently. Therefore, the likelihood of the observed values in $\mathbf{Y}$ is the product of $d$ independent GPs:
\[\text{p}(\mathbf{Y}| \mathbf{X}, \boldsymbol{\theta})= \prod_{i=1}^d \frac{\exp\left[-\frac{1}{2} \mathbf{y}_{:,i}^\text{T} \mathbf{K}^{-1} \mathbf{y}_{:,i} \right]} {(2\pi)^{n/2}|\mathbf{K}|^{1/2}}.\]
Then we can adjust $\mathbf{X}$ and $\boldsymbol{\theta}$ to maximize this likelihood. If we were only adjusting $\mathbf{X}$, it would be akin to determining the shadow which is most likely to correspond to your body's shape. Because we're adjusting $\boldsymbol{\theta}$ as well, imagine manipulating the light source and the surface for the shadow to appear on as well.

\begin{figure}[h]
{\includegraphics[width=12cm]{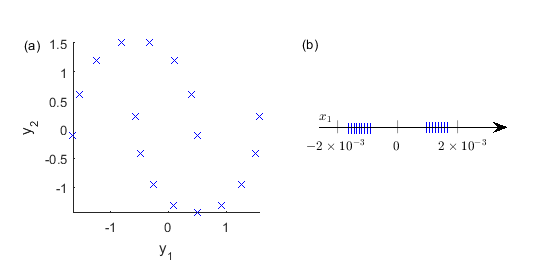}}  
{\caption{\textbf{(a)} Seventeen data points, $\mathbf{Y}$, each of two dimensions. \newline 
\color{white}{d }\color{black}  \textbf{(b)} A one-dimensional projection, $\mathbf{X}$, of 
those points.}   
\label{fig:toyDR} }
\end{figure}

\begin{figure}[b]
 \begin{center}
  \includegraphics[height=8.9cm] {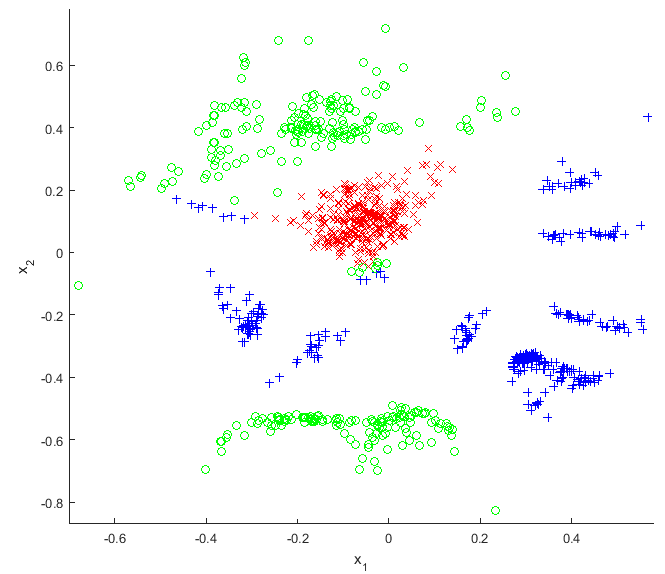}
  \caption{Although the GP-LVM algorithm uses unsupervised learning, the three classes of data in this originally 12-D dataset generally occupy different regions.}
  \label{fig:oil}
 \end{center}
\end{figure}

Figure~\ref{fig:toyDR} gives an example outcome of the most likely $\mathbf{X}$ for a certain $\mathbf{Y}$, with $n=17$, $d=2$, and $q=1$. $\mathbf{Y}$ was preprocessed to have zero mean and the same variance in each dimension. Using an optimizer based on scaled conjugate gradients, the optimal $\boldsymbol{\theta}$ and $\mathbf{X}$ were then found. The former is $\{\sigma,l,\beta\} = \{1.05,3.3\times 10^{-4},93\}$ and the latter is given in Figure~\ref{fig:toyDR}(b). In this example, $\mathbf{X}$ is more than just a `shadow' (simple projection) of $\mathbf{Y}$: the above method of dimensionality reduction is powerful enough to have learned that the points in the left curve in Figure~\ref{fig:toyDR}(a) should be grouped together while the points in the right curve form a separate group; no simple lamplight projection can achieve this. The price we pay for this flexibility is in the high number of parameters being fit: the total here is $3+17 \times 1 = 20$, i.e.\ $\boldsymbol{\theta}$ plus the one-dimensional values in $\mathbf{X}$. Usually $\mathbf{X}$ are referred to as latent variables rather than parameters, and our overall setup is called the `Gaussian process latent variable model' (GP-LVM). The technique was first presented by Neil Lawrence in 2004; he examined a set of oil-flow data in which $n=1000$ and $d=12$, which we reduce to $q=2$ in Figure~\ref{fig:oil}.

This tutorial on the GP-LVM is provided by Mind Foundry, a technology spin-out from the University of Oxford. 
%Computer code: \textit{github.com/mebden/GPLVMtutorial}

\captionsetup[figure]{justification=centering}

% GLtutorialReport.m, GLtuteData.mat
%More general GP-LVM code is available at \textit{http://sheffieldml.github.io/software.html}.

\end{document}